\documentclass[12pt]{amsart} 
\usepackage{amssymb,amsmath,amscd} 
\def\mysavedown#1{\edef\mysubs{\mysubs#1}}
\def\mysaveup#1{\edef\mysups{\mysups#1}}
\def\mydown#1{{\mytensor}_{\vphantom{\mysubs}#1}}
\def\myup#1{{\mytensor}^{\vphantom{\mysups}#1}}
\def\tensor#1#2{
  #1
  \def\mytensor{\vphantom{#1}}
  \def\mysubs{\relax}
  \def\mysups{\relax}
  \let\down=\mysavedown
  \let\up=\mysaveup
  #2
  \let\down=\mydown
  \let\up=\myup
  #2
  }

\theoremstyle{plain} 
\newtheorem{Lem}{Lemma}[section] 
\newtheorem{Prop}[Lem]{Proposition} 
\newtheorem{Thm}[Lem]{Theorem} 
 
\theoremstyle{definition} 
\newtheorem{Def}[Lem]{Definition} 
\newtheorem{Rem}[Lem]{Remark} 
 
\newtheorem{Prbl}[Lem]{Problem} 
\newtheorem{Proc}[Lem]{Procedure} 
\errorcontextlines=0  \numberwithin{equation}{section}

\newcommand{\bbC}{{\mathbb C}} 
\newcommand{\bbR}{{\mathbb R}}

\newcommand{\bbN}{{\mathbb N}} 
\newcommand{\bbK}{{\mathbb K}} 
\newcommand{\calA}{{\mathcal A}} 
\newcommand{\calB}{{\mathcal B}} 
\newcommand{\calH}{{\mathcal H}} 
\newcommand{\calL}{{\mathcal L}}

\newcommand{\calP}{{\mathcal P}} 
\newcommand{\calS}{{\mathcal S}} 
\newcommand{\calT}{{\mathcal T}} 
\newcommand{\calV}{{\mathcal V}} 
 
\newcommand{\id}{\mathrm{id}} 
\newcommand{\Span}{\mathrm{Span}} 
\newcommand{\frakl}{{\mathfrak l}} 
\newcommand{\frakr}{{\mathfrak r}} 

\newcommand{\fR}{{\mathfrak R}}
\newcommand{\fS}{{\mathfrak S}}
 
\newcommand{\bpi}{\begin{picture}} 
\newcommand{\epi}{\end{picture}}

\begin{document} 

\title[(2,1)-generators of algebraic curvature tensors] 
{Generators of algebraic curvature tensors based on a (2\,1)-symmetry} 
\author[B. Fiedler]{Bernd Fiedler}
\address{Bernd Fiedler \\ Mathematisches Institut \\ Universit\"at Leipzig\\ 
Augustusplatz 10/11 \\ D-04109 Leipzig \\ Germany}
\urladdr{http://home.t-online.de/home/Bernd.Fiedler.RoschStr.Leipzig/}  
\email{Bernd.Fiedler.RoschStr.Leipzig@t-online.de}
\subjclass{53B20, 15A72, 05E10, 16D60, 05-04} 

\begin{abstract}
We consider generators of algebraic curvature tensors ${\mathfrak R}$ which can be constructed by a Young symmetrization of product tensors $U\otimes w$ or $w\otimes U$, where $U$ and $w$ are covariant tensors of order 3 and 1. We assume that $U$ belongs to a class of the infinite set ${\mathfrak S}$ of irreducible symmetry classes characterized by the partition $(2\,1)$. We show that the set ${\mathfrak S}$ contains exactly one symmetry class $S_0\in{\mathfrak S}$ whose elements $U\in S_0$ can not play the role of generators of tensors ${\mathfrak R}$. The tensors $U$ of all other symmetry classes from ${\mathfrak S}\setminus\{S_0\}$ can be used as generators for tensors ${\mathfrak R}$.

Using Computer Algebra we search for such generators whose coordinate representations are polynomials with a minimal number of summands. For a generic choice of the symmetry class of $U$ we obtain lengths of 8 summands. In special cases these numbers can be reduced to the minimum 4. If this minimum occurs then $U$ admits an index commutation symmetry. Furthermore minimal lengths are possible if $U$ is formed from torsion-free covariant derivatives of alternating 2-tensor fields.

We apply ideals and idempotents of group rings ${\mathbb C}[{\calS}_r]$ of symmetric groups $\calS_r$, Young symmetrizers, discrete Fourier transforms and Littlewood-Richardson products. For symbolic calculations we used the Mathematica packages {\sf Ricci} and {\sf PERMS}.
\end{abstract}

\maketitle 

%
%

\section{Introduction}
In \cite{fie03b,fie03c} we constructed and investigated generators of \textit{algebraic covariant derivative curvature tensors} which contained 3-times covariant tensors $U$ with a so-called $(2\,1)$-symmetry. Later we discovered that the same constructions can be applied in the simpler case of \textit{algebraic curvature tensors}, too. In the present paper we carry out all constructions and investigations of \cite{fie03b,fie03c} for algebraic curvature tensors and join the results about these tensors to the results of \cite{fie03b,fie03c}.

Let ${\calT}_r V$ be the vector space of the $r$-times covariant tensors $T$ over a finite-dimensional $\bbK$-vector space $V$, $\bbK = \bbR$ or  $\bbK = \bbC$. We assume that $V$ possesses a {\itshape fundamental tensor} $g \in {\calT}_2 V$ (of arbitrary signature) which can be used for raising and lowering of tensor indices.
\begin{Def}
Let $\calA(V)\subset\calT_4 V$ and $\calA'(V)\subset\calT_5 V$ be the spaces of all \textit{algebraic curvature tensors} and all \textit{algebraic covariant derivative curvature tensors}, respectively, i.e. those tensors
$\fR\in\calA(V)$, $\fR'\in\calA'(V)$ which satisfy
\begin{eqnarray}
 & & \fR(w,x,y,z) \;=\; - \fR(w,x,z,y) \;=\; \fR(y,z,w,x)\,, \\
 & & \fR(w,x,y,z) + \fR(w,y,z,x) + \fR(w,z,x,y) \;=\; 0\,, \\
 & & \fR'(w,x,y,z;u) \;=\; - \fR'(w,x,z,y;u) \;=\; \fR'(y,z,w,x;u)\,, \\
 & & \fR'(w,x,y,z;u) + \fR'(w,y,z,x;u) + \fR'(w,z,x,y;u) \;=\; 0\,, \\
 & & \fR'(w,x,y,z;u) + \fR'(w,x,z,u;y) + \fR'(w,x,u,y;z) \;=\; 0
\end{eqnarray}
for all $u, w, x, y, z\in V$.
\end{Def}
$\fR$ and $\fR'$ have the symmetries of the Riemann tensor $R$ of a Levi-Civita connection $\nabla$ and the covariant derivative $\nabla R$.
%
%
%

Let $\calS^p(V)$, ${\Lambda}^p(V)$ be the spaces of totally symmetric/alternating $p$-forms over $V$. Then the following tensors
\begin{eqnarray}
\gamma (S)_{i j k l} & := &
S_{i l} S_{j k} - S_{i k} S_{j l}\,,\\
\alpha (A)_{i j k l} & := &
2\,A_{i j} A_{k l} + A_{i k} A_{j l}
-  A_{i l} A_{j k}\,, \\
\hspace*{1cm}\hat{\gamma} (S,\hat{S})_{i j k l s} & := &
S_{il}{\hat{S}}_{jks} - S_{jl}{\hat{S}}_{iks} + S_{jk}{\hat{S}}_{ils} - S_{ik}{\hat{S}}_{jls}\,, \label{e1.8}\\
 & & S\in\calS^2(V)\,,\,\hat{S}\in\calS^3(V)\,,\,A\in {\Lambda}^2(V)\,, \nonumber
\end{eqnarray}
are generators of $\calA(V)$, $\calA'(V)$. P. Gilkey \cite[pp.41--44, p.236]{gilkey5} and B. Fiedler \cite{fie20,fie03b} gave different proofs for
\begin{Thm} $\,$\label{thm1.2}%
\begin{enumerate}
\item{$\calA(V) = \Span_{S\in\calS^2(V)}\{\gamma(S)\} = \Span_{A\in {\Lambda}^2(V)}\{\alpha(A)\}$.}
\item{$\calA'(V) = \Span_{S\in\calS^2(V) , \hat{S}\in\calS^3(V)}\{\hat{\gamma}(S , \hat{S})\}$. \label{item1.2.2}}
\end{enumerate}
\end{Thm}
The tensors $\gamma(S)$, $\alpha(A)$ and $\hat{\gamma}(S , \hat{S})$ are expressions which arise from
$S\otimes S$, $A\otimes A$, $S\otimes \hat{S}$ or $\hat{S}\otimes S$ by a symmetrization
\begin{eqnarray}
 & & \gamma(S)\;=\;\textstyle{\frac{1}{12}}\,y_t^{\ast}(S\otimes S)
\;\;\;,\;\;\;
\alpha(A)\;=\;\textstyle{\frac{1}{12}}\,y_t^{\ast}(A\otimes A) \\
 & & \hat{\gamma}(S,\hat{S})\;=\;\textstyle{\frac{1}{4}}\,y_{t'}^{\ast}(S\otimes \hat{S})\;=\;\textstyle{\frac{1}{4}}\,y_{t'}^{\ast}(\hat{S}\otimes S)
\end{eqnarray}
where $y_t$, $y_{t'}$ are the \textit{Young symmetrizers} of the \textit{Young tableaux}
\begin{eqnarray} \label{e1.11}%
t\;= \;
\begin{array}{|c|c|}
\hline
1 & 3 \\
\hline
2 & 4 \\
\hline
\end{array}
& \;\;\;,\;\;\; &
t'\;=\;
\begin{array}{|c|c|c|c}
\cline{1-3}
1 & 3 & 5 & \\
\cline{1-3}
2 & 4 &\multicolumn{2}{c}{\;\;\;} \\
\cline{1-2}
\end{array}
\end{eqnarray}
(See \cite{fie20,fie03b}. See also Section \ref{sec2} for definitions.).

In the present paper we search for similar generators of $\calA(V)$ and $\calA'(V)$ which, however, are based on other product tensors. We use Boerner's definition of {\itshape symmetry classes} for tensors $T\in\calT_r V$ by right ideals $\frakr\subseteq\bbK[\calS_r]$ of the group ring $\bbK[\calS_r]$ of the symmetric group $\calS_r$ (see Section \ref{sec2} and \cite{boerner,boerner2,fie16,fie18}). On this basis we investigate the following
\begin{Prbl} \label{probl1.3}%
We search for generators of $\calA(V)$, $\calA'(V)$ which can be formed by a suitable symmetry operator from tensors
\begin{eqnarray} \label{e1.12}%
\calA(V)\,:\;\;\;U\otimes w\;\;,\;\;w\otimes U\hspace{0.3cm} &\;\;\;,\;\;\;& U\in\calT_3 V\;,\;w\in\calT_1 V\,,\\
\calA'(V)\,:\;\;\;U\otimes W\;\;,\;\;W\otimes U &\;\;\;,\;\;\;& U\in\calT_3 V\;,\;W\in\calT_2 V\,, \label{e1.13}%
\end{eqnarray}
where $W$ and $U$ belong to symmetry classes of $\calT_2 V$ and $\calT_3 V$ which are defined by minimal right ideals $\frakr\subset\bbK[\calS_2]$ and $\hat{\frakr}\subset\bbK[\calS_3]$, respectively.
\end{Prbl}
Here is a summary of our main results. The subject of the present paper is the determination of generators (\ref{e1.12}) of $\calA(V)$. However, for comparing purposes we repeat also results concerning generators (\ref{e1.13}) of $\calA'(V)$ which were proved in \cite{fie03b,fie03c}.
\begin{Thm} \label{thm1.4}%
A solution of Problem {\rm\ref{probl1.3}} can be constructed at most from such products {\rm (\ref{e1.12})}or {\rm (\ref{e1.13})} whose factors $U\in\calT_3 V$, $W\in\calT_2 V$ lie in symmetry classes which belong\footnote{See Section \ref{sec2} and \cite{boerner,boerner2,fie03b} for the connection between partitions and symmetry classes or right ideals respectively.} to the following partitions of {\rm 3} or {\rm 2}:
\begin{center}
{\rm
\begin{tabular}{|l|c|l|l|}
\hline
product & & partitions for $U$, $W$ & \\
\hline
$\fR:\,U\otimes w$ & (a) & $U\leftrightarrow\,(2\,1)$ & \\
\hspace*{0.8cm}$w\otimes U$ & & & \\
\hline
$\fR':\,U\otimes W$ & (a') & $U\leftrightarrow\,(3)\hspace{0.3cm}\;\mathrm{and}\;W\leftrightarrow\,(2)$ & $U$ and $W$ symmetric \\
\hspace*{0.85cm}$W\otimes U$ & (b') & $U\leftrightarrow\,(2\,1)\;\mathrm{and}\;W\leftrightarrow\,(2)$ & $W$ symmetric \\
 & (c') & $U\leftrightarrow\,(2\,1)\;\mathrm{and}\;W\leftrightarrow\,(1^2)$ & $W$ skew-symmetric \\
\hline
\end{tabular}\,.
}
\end{center}
\end{Thm}
The case (a') of Theorem \ref{thm1.4} is specified by Theorem \ref{thm1.2}(\ref{item1.2.2}) and (\ref{e1.8}) (see \cite{fie03b}).
The cases (a), (b') and (c') of Theorem \ref{thm1.4} lead to
\begin{Thm} \label{thm1.5}%
Let $\frakr\subset\bbK[\calS_3]$ be a minimal right ideal belonging to the partition $(2\,1)\vdash 3$ and let $\calT_{\frakr}$ be the symmetry class of tensors $U\in\calT_3 V$ that is defined by $\frakr$. Then the following statements are equivalent:
\begin{enumerate}
\item{$\calA(V) = \Span_{U\in\calT_{\frakr}\,,\,w\in\calT_1 V}\{y_t^{\ast}(U\otimes w)\}$. \label{i1.5.1}}
\item{$\calA'(V) = \Span_{U\in\calT_{\frakr}\,,\,S\in\calS^2(V)}\{y_{t'}^{\ast}(U\otimes S)\}$. \label{i1.5.2}}
\item{$\calA'(V) = \Span_{U\in\calT_{\frakr}\,,\,A\in{\Lambda}^2(V)}\{y_{t'}^{\ast}(U\otimes A)\}$. \label{i1.5.3}}
\item{The right ideal $\frakr$ is different from the right ideal
$\frakr_0 = f_0\cdot\bbK[\calS_3]$ which is generated by the idempotent \label{i1.5.4}
\begin{eqnarray} \label{e1.14}%
f_0 & := & \textstyle{\frac{1}{2}}\,\{\id - (1\,3)\} - \textstyle{\frac{1}{6}}\,y \;\;\;\;,\;\;\;\;
y\;:=\;\sum_{p\in\calS_3}\,\mathrm{sign}(p)\,p\,.
\end{eqnarray}
}
\end{enumerate}
Here $y_t$ and $y_{t'}$ are the Young symmetrizers of the Young tableaux {\rm (\ref{e1.11})}.
\end{Thm}
The Statements (\ref{i1.5.1}), (\ref{i1.5.2}), (\ref{i1.5.3}) of Theorem \ref{thm1.5} are independent of the order of the factors $U, S, A, w$ since the following Lemma holds.
\begin{Lem} \label{lem1.6}%
Let $\calT_{\frakr}$ be the symmetry class and let $y_t\,,\,y_{t'}$ be the Young symmetrizers from Theorem {\rm \ref{thm1.5}}.
If $U\in\calT_{\frakr}\,,\,S\in\calS^2(V)\,,\,A\in{\Lambda}^2(V)$ and $w\in\calT_1 V$, then it holds
\begin{eqnarray}
y_t^{\ast}(U\otimes w) & = & y_t^{\ast}(w\otimes U)\,, \label{e1.15}\\
y_{t'}^{\ast}(U\otimes S) & = & y_{t'}^{\ast}(S\otimes U)\,, \label{e1.16}\\
y_{t'}^{\ast}(U\otimes A) & = & - y_{t'}^{\ast}(A\otimes U)\,. \label{e1.17}
\end{eqnarray}
\end{Lem}
\begin{Rem}
The set $\fS$ of symmetry classes $\calT_{\frakr}$ considered in Theorem \ref{thm1.5} is an infinite set. Theorem \ref{thm1.5} says that exactly the tensors $U$ from symmetry classes $\calT_{\frakr}\in\fS\setminus\{\calT_{\frakr_0}\}$ yield generators of $\calA(V)$ or $\calA'(V)$ respectively.

The equivalence of the statements (\ref{i1.5.2}), (\ref{i1.5.3}), (\ref{i1.5.4}) of Theorem \ref{thm1.5} was shown already in \cite{fie03b}. It is remarkable that we also have to exclude a single symmetry class from $\fS$ if we search for generators of $\calA(V)$, and that the forbidden symmetry class is the same class $\calT_{\frakr_0}$ which had to be excluded in the construction of generators of $\calA'(V)$.
\end{Rem}
For the generators (\ref{i1.5.1}), (\ref{i1.5.2}), (\ref{i1.5.3}) of Theorem \ref{thm1.5} we can also determine operators of the type $\alpha$, $\gamma$, $\hat{\gamma}$ which yield the coordinate representation of these generators.
 However, these operators depend now on the right ideal $\frakr$ (or its generating idempotent $e$) that defines the symmetry class of $U$. And they yield no short expressions of 2, 3, or 4 terms but longer expressions. The {\it search for shortest expressions of this type} is a further subject of our paper. Some of our main results are collected in
\begin{Thm} \label{thm1.8}%
Consider the situation of Theorem {\rm\ref{thm1.5}} and assume $\dim V\ge 3$, $\frakr\not=\frakr_0$. Then it holds:
\begin{enumerate}
\item{The coordinates of $y_t^{\ast}(U\otimes w)$, $y_{t'}^{\ast}(U\otimes S)$, $y_{t'}^{\ast}(U\otimes A)$ are sums of the following lengths
\begin{center}
\begin{tabular}{|c|l|c|c|c|}
\hline
 & & $y_t^{\ast}(U\otimes w)$ & $y_{t'}^{\ast}(U\otimes S)$ & $y_{t'}^{\ast}(U\otimes A)$ \\
\hline
{\rm (a)} & {\rm generic case for $\frakr$} & {\rm 8} & {\rm 16} & {\rm 20} \\
\hline
{\rm (b)} & {\rm $\frakr$ producing minimal length} & {\rm 4} & {\rm 12} & {\rm 10} \\
\hline
\end{tabular}
\end{center}
}
\item{There exist exactly {\rm 2} minimal right ideals $\frakr\not=\frakr_0$ of $(2\,1)\vdash 3$ which lead to the minimal lengths of case (b) for all tensors
$y_t^{\ast}(U\otimes w)$, $y_{t'}^{\ast}(U\otimes S)$ and $y_{t'}^{\ast}(U\otimes A)$.}
\item{If the coordinates of $y_t^{\ast}(U\otimes w)$, $y_{t'}^{\ast}(U\otimes S)$, $y_{t'}^{\ast}(U\otimes A)$ have the minimal lengths of case {\rm (b)} then $U$ admits an index commutation symmetry.}
\end{enumerate}
\end{Thm}
Further results are given in Section \ref{sec4}.
\begin{Rem}
The concept ''expression of minimal length'' depends on the method which we use to reduce expressions (see Section \ref{sec4}, Procedure \ref{proc4.3}). In \cite[Remark 3.9]{fie03c} we discuss a generalization of our reduction method which could possibly lead to a further decrease of the numbers in Theorem \ref{thm1.8}.
\end{Rem}
Examples of tensors $U$ with a $(2\,1)$-symmetry from ${\mathfrak S}$ can be constructed from covariant derivatives of certain tensor fields. In \cite{fie03a,fie03c} we proved
\begin{Prop} {\rm\bf (Examples of $(2\,1)$-symmetries)}\\
Let $\nabla$ be a torsion-free covariant derivative on a $C^{\infty}$-manifold $M$, $\dim M\ge 2$. Further let $\psi\in\calS^2 M$, $\omega\in{\Lambda}^2 M$ be differentiable tensor fields of order {\rm 2} on $M$ which are symmetric or skew-symmetric, respectively. Then the infinite set $\fS$ of symmetry classes contains {\rm 2} classes $\calT_{\frakr_s} , \calT_{\frakr_a}\in\fS$ such that
\begin{eqnarray} \label{e1.18}%
\hspace*{1cm}\forall\,p\in M:\;\;(\nabla\psi - \mathrm{sym}(\nabla\psi))|_p\in\calT_{\frakr_s}
& \;\;\;,\;\;\; &
(\nabla\omega - \mathrm{alt}(\nabla\omega))|_p\in\calT_{\frakr_a}\,.
\end{eqnarray}
Here '{\rm sym}' denotes the symmetrization and '{\rm alt}' the anti-symmetrization of tensors.
\end{Prop}
More details can be found in Remark \ref{rem2.12}. For tensors (\ref{e1.18}) we obtained
\begin{Thm} \label{thm1.10}%
If we consider tensors $U\in\calT_{\frakr}$, $S\in\calS^2(V)$, $A\in{\Lambda}^2(V)$, $w\in\calT_1 V$ on a tangent space $V = T_p M$ of a differentiable manifold $M$, $\dim M\ge 3$, and form $U$ by one of the formulas {\rm (\ref{e1.18})}, then we obtain the shortest lengths from Theorem {\rm \ref{thm1.8}, (1b)} exactly in the following cases:
\begin{enumerate}
\item{$y_{t'}^{\ast}(U\otimes S)$ and $U = (\nabla\psi - {\rm sym}(\nabla\psi))|_p$, $\psi\in\calT_2 M$ symmetric,}
\item{$y_{t}^{\ast}(U\otimes w)$, $y_{t'}^{\ast}(U\otimes S)$, $y_{t'}^{\ast}(U\otimes A)$ and
$U = (\nabla\omega - {\rm alt}(\nabla\omega))|_p$, $\omega\in\calT_2 M$ skew-symmetric.}
\end{enumerate}
\end{Thm}
Here is a brief outline to the paper. In Section \ref{sec2} we give a summary of basic facts about symmetry classes, Young symmetrizers and discrete Fourier transforms. These tools are needed to obtain the infinite set ${\mathfrak S}$ of symmetry classes for $U$. In Section \ref{sec3} we prove the Theorems \ref{thm1.4}, \ref{thm1.5} and Lemma \ref{lem1.6} using the \textit{Littlewood-Richardson rule} and computer calculations with group ring elements and tensor coordinates. In Section \ref{sec4} we construct short coordinate representations for the tensors $y_t^{\ast}(U\otimes w)$ by determining and solving a complete system of linear identities for the tensors $U$. Furthermore we prove the Theorems \ref{thm1.8} and \ref{thm1.10} in this Section.

Many results were obtained by computer calculations by means of the {\sf Mathematica} packages {\sf Ricci} \cite{ricci3} and {\sf PERMS} \cite{fie10}. The {\sf Mathematica} notebooks of these calculations are available at \cite{fie21}.\vspace{10pt}

\section{Symmetry classes, Young symmetricers, discrete Fourier transforms} \label{sec2}
The vector spaces $\calA(V)$ and $\calA'(V)$ 
are {\itshape symmetry classes} in the sence of H. Boerner \cite[p.127]{boerner}. We denote by $\bbK [{\calS}_r]$ the {\itshape group ring} of a symmetric group ${\calS}_r$ over the field $\bbK$. Every group ring element $a = \sum_{p \in {\calS}_r} a(p)\,p \in \bbK [{\calS}_r]$ acts as so-called {\itshape symmetry operator} on tensors $T \in {\calT}_r V$ according to the definition
\begin{eqnarray}
(a T)(v_1 , \ldots , v_r) & := & \sum_{p \in {\calS}_r} a(p)\,
T(v_{p(1)}, \ldots , v_{p(r)}) \;\;\;\;\;,\;\;\;\;\;
v_i \in V \,. \label{e2.1}
\end{eqnarray}
Equation \eqref{e2.1} is equivalent to
$(a T)_{i_1 \ldots i_r} = \sum_{p \in {\calS}_r} a(p)\,
T_{i_{p(1)} \ldots  i_{p(r)}}$.
\begin{Def}
Let $\frakr \subseteq \bbK [{\calS}_r]$ be a right ideal of $\bbK [{\calS}_r]$ for which an $a \in \frakr$ and a $T \in {\calT}_r V$ exist such that $aT \not= 0$. Then the tensor set
\begin{eqnarray}
{\calT}_{\frakr} & := & \{ a T \;|\; a \in \frakr \;,\;
T \in {\calT}_r V \}
\end{eqnarray}
is called the {\itshape symmetry class} of tensors defined by $\frakr$.
\end{Def}
Since $\bbK [{\calS}_r]$ is semisimple for $\bbK = \bbR , \bbC$, every right ideal $\frakr \subseteq \bbK [{\calS}_r]$ possesses a generating idempotent $e$, i.e. $\frakr$ fulfils $\frakr = e \cdot \bbK [{\calS}_r]$. It holds (see e.g. \cite{fie20} or \cite{boerner,boerner2})
\begin{Lem}
If $e$ is a generating idempotent of $\frakr$, then a tensor $T \in {\calT}_r V$ belongs to ${\calT}_{\frakr}$ iff
$e T = T$.
Thus we have
${\calT}_{\frakr} = \{ eT \;|\; T \in {\calT}_r V \}$.
\end{Lem}
Now we summarize tools from our Habilitationsschrift \cite{fie16} (see also its summary \cite{fie18}).
We make use of the following connection
between $r$-times covariant tensors $T \in {\calT}_r V$
and elements of the {\it group ring}
${\bbK} [{\calS}_r]$.
\begin{Def} \label{def2.3}
 Any tensor
 $T \in {\calT}_r V$
 and any $r$-tuple
 $b := (v_1 , \ldots , v_r ) \in V^r$
 of
 $r$
 vectors from
 $V$
 induce a function
 $T_b : {\calS}_r \rightarrow {\bbK}$
 according to the rule
 \begin{eqnarray}
T_b (p) & := & T(v_{p(1)} , \ldots , v_{p(r)})\;\;\;,\;\;\;p \in {\calS}_r \,.
\end{eqnarray}
We identify this function with the group ring element
$T_b := \sum_{p \in {\calS}_r}T_b (p)\,p \in {\bbK} [{\calS}_r]$.
\end{Def}
Obviously,
two tensors $S , T \in {\calT}_r V$ fulfil $S = T$ iff
$S_b = T_b$ for all
$b \in V^r$.
We denote by '$\ast$'
the mapping
$\ast : a = \sum_{p \in {\calS}_r} a(p)\,p \;\mapsto\; a^{\ast} :=
\sum_{p \in {\calS}_r} a(p)\,p^{-1}$. Then the following important formula\footnote{See B. Fiedler \cite[Sec.III.1]{fie16} and B. Fiedler \cite{fie17}.} holds
\begin{eqnarray} \label{e2.4}
\forall\,T\in\calT_r V\;,\;a\in\bbK[\calS_r]\;,\;b\in V^r\;:\;\;\;\;
(a\,T)_b & = & T_b\cdot a^{\ast}\,.
\end{eqnarray}
Now it can be shown that all $T_b$ of tensors $T$ of a given symmetry class lie in a certain left ideal of ${\bbK}[{\calS}_r]$.
\begin{Prop}\hspace{-1mm}\footnote{See B. Fiedler \cite{fie17} or
B. Fiedler \cite[Prop. III.2.5, III.3.1, III.3.4]{fie16}.}
\label{prop2.4}%
Let $e \in {\bbK}[{\calS}_r]$ be an idempotent. Then a
$T \in {\calT}_r V$ 
fulfils the condition
$eT = T$
iff
$T_b \;\in\; {\frakl} := {\bbK} [{\calS}_r] \cdot e^{\ast}$ for all
$b \in V^r$, i.e.
all $T_b$ of $T$
lie in the left ideal ${\frakl}$ generated by $e^{\ast}$.
\end{Prop}
The proof follows easily from (\ref{e2.4}). Since a rigth ideal $\frakr$ defining a symmetry class and the left ideal $\frakl$ from Proposition \ref{prop2.4} satisfy $\frakr = \frakl^{\ast}$, we denote symmetry classes also by $\calT_{\frakl^{\ast}}$. A further result is
\begin{Prop}\hspace{-1mm}\footnote{See B. Fiedler \cite{fie17} or
B. Fiedler \cite[Prop. III.2.6]{fie16}.}
\label{prop2.5}%
If $\dim V \ge r$, then every left ideal
${\frakl} \subseteq {\bbK}[{\calS}_r]$ fulfils
${\frakl} = {\calL}_{\bbK} \{ T_b \;|\;
T \in {\calT}_{{\frakl}^{\ast}} \,,\, b \in V^r \}$.
(Here ${\calL}_{\bbK}$ denotes the forming of the linear closure.)
\end{Prop}
If $\dim V < r$, then the $T_b$ of the tensors from
${\calT}_{{\frakl}^{\ast}}$ 
will span only a linear subspace of
${\frakl}$ 
in general.

Important special symmetry operators are Young symmetrizers, which are defined by means of Young tableaux.

A {\itshape Young tableau} $t$ of $r\in\bbN$ is an arrangement of $r$ boxes such that
\begin{enumerate}
\item{the numbers ${\lambda}_i$ of boxes in the rows $i = 1 , \ldots , l$ form a decreasing sequence
${\lambda}_1 \ge {\lambda}_2 \ge \ldots \ge {\lambda}_l > 0$ with
${\lambda}_1 + \ldots + {\lambda}_l = r$,}
\item{the boxes are fulfilled by the numbers $1, 2, \ldots , r$ in any order.}
\end{enumerate}
For instance, the following graphics shows a Young tableau of $r = 16$.
\[\left.
\begin{array}{cc|c|c|c|c|c|c}
\cline{3-7}
{\lambda}_1 = 5 & \;\;\; & 11 & 2 & 5 & 4 & 12 & \\
\cline{3-7}
{\lambda}_2 = 4 & \;\;\; & 9 & 6 & 16 & 15 & \multicolumn{2}{c}{\;\;\;} \\
\cline{3-6}
{\lambda}_3 = 4 & \;\;\; & 8 & 14 & 1 & 7 & \multicolumn{2}{c}{\;\;\;} \\
\cline{3-6}
{\lambda}_4 = 2 & \;\;\; & 13 & 3 & \multicolumn{4}{c}{\hspace{2cm}} \\
\cline{3-4}
{\lambda}_5 = 1 & \;\;\; & 10 & \multicolumn{4}{c}{\hspace{2cm}} \\
\cline{3-3}
\end{array}\right\}\;=\;t\,.
\]
Obviously, the unfilled arrangement of boxes, the {\itshape Young frame}, is characterized by a partition
$\lambda = ({\lambda}_1 , \ldots , {\lambda}_l) \vdash r$ of $r$.

If a Young tableau $t$ of a partition $\lambda \vdash r$ is given, then the {\itshape Young symmetrizer} $y_t$ of $t$ is defined by\footnote{We use the convention $(p \circ q) (i) := p(q(i))$ for the product of two permutations $p, q$.}
\begin{eqnarray}
y_t & := & \sum_{p \in {\calH}_t} \sum_{q \in {\calV}_t} \mathrm{sign}(q)\, p \circ q
\end{eqnarray}
where ${\calH}_t$, ${\calV}_t$ are the groups of the {\itshape horizontal} or
{\itshape vertical permutations} of $t$ which only permute numbers within rows or columns of $t$, respectively. The Young symmetrizers of $\bbK [{\calS}_r]$ are essentially idempotent and define decompositions
\begin{eqnarray}
\bbK [{\calS}_r] \;=\;
\bigoplus_{\lambda \vdash r} \bigoplus_{t \in {\calS\calT}_{\lambda}}
\bbK [{\calS}_r]\cdot y_t
& \;\;,\;\; &
\bbK [{\calS}_r] \;=\;
\bigoplus_{\lambda \vdash r} \bigoplus_{t \in {\calS\calT}_{\lambda}}
y_t \cdot \bbK [{\calS}_r] \label{e2.6}
\end{eqnarray}
of $\bbK [{\calS}_r]$ into minimal left or right ideals. In \eqref{e2.6}, the symbol ${\calS\calT}_{\lambda}$ denotes the set of all standard tableaux of the partition $\lambda$. \textit{Standard tableaux} are Young tableaux in which the entries of every row and every column form an increasing number sequence.\footnote{About Young symmetrizers and
Young tableaux see for instance
\cite{boerner,boerner2,full4,fulton,jameskerb,kerber,littlew1,mcdonald,%
muell,naimark,%
waerden,weyl1}. In particular, properties of Young symmetrizers in the case
${\bbK} \not= {\bbC}$ are described in \cite{muell}.}

Every \textit{irreducible character} $\chi: \calS_r\rightarrow\bbC$ of $\calS_r$ induces a \textit{centrally primitive idempotent}
$\chi := \frac{\chi({\rm id})}{r!}\sum_{p\in\calS_r}\,\chi(p)\,p$ which generates a \textit{minimal two-sided ideal} ${\mathfrak a}_{\chi} := \bbK[\calS_r]\cdot\chi$. There is a \textit{one-to-one correspondence} $\chi\leftrightarrow\lambda$ between the $\chi$ and the partitions $\lambda\vdash r$.
For every $\chi$ there exists a unique $\lambda\vdash r$ such that
\begin{eqnarray}
{\mathfrak a}_{\chi} & = &
\bigoplus_{t \in {\calS\calT}_{\lambda}}
\bbK [{\calS}_r]\cdot y_t
\;=\;
\bigoplus_{t \in {\calS\calT}_{\lambda}}
y_t \cdot \bbK [{\calS}_r]
\end{eqnarray}
The set of all Young symmetrizers $y_t$ which lie in ${\mathfrak a}_{\chi}$ is equal to the set of all $y_t$ whose tableau $t$ has a frame $\lambda\vdash r$.
Furthermore two minimal left ideals $\frakl_1, \frakl_2\subseteq\bbK[\calS_r]$ or two minimal right ideals $\frakr_1, \frakr_2\subseteq\bbK[\calS_r]$ are \textit{equivalent} iff they lie in the same ideal ${\mathfrak a}_{\chi}$.
Now we say that a symmetry class $\calT_{\frakr}$ \textit{belongs to} $\lambda\vdash r$ iff
$\frakr\subseteq{\mathfrak a}_{\chi}$ and $\chi$ corresponds to $\lambda$. Then we write also ${\mathfrak a}_{\lambda}$ for ${\mathfrak a}_{\chi}$.

S.A. Fulling, R.C. King, B.G.Wybourne and C.J. Cummins showed in \cite{full4} that the symmetry classes of the Riemannian curvature tensor $R$ and its {\itshape symmetrized\footnote{$(\,\ldots\,)$ denotes the symmetrization with respect to the indices $s_1, \ldots , s_u$.} covariant derivatives}
\begin{eqnarray}
\left({\nabla}^{(u)}R\right)_{i j k l s_1 \ldots s_u} & := & {\nabla}_{(s_1} {\nabla}_{s_2} \ldots {\nabla}_{s_u)} R_{i j k l}\;=\;R_{i j k l\,;\,(s_1 \ldots s_u)}
\end{eqnarray}
are generated by special Young symmetrizers\footnote{A proof of this result of \cite{full4} can be found in \cite[Sec.6]{fie12}, too. See also \cite{fie03c} for more details.}. 
In the present paper we use only 
\begin{Thm} {\rm\bf (Fulling, King, Wybourne, Cummins)} \label{thm2.6} \\
Let $y_t$, $y_{t'}$ be the Young symmetrizers of the standard tableaux
\begin{eqnarray} \label{e2.8}%
t\;=\;
\begin{array}{|c|c|}
\hline
1 & 3 \\
\hline
2 & 4 \\
\hline
\end{array}
& \;\;\;,\;\;\;
t'\;=\;
\begin{array}{|c|c|c|c}
\cline{1-3}
1 & 3 & 5 & \\
\cline{1-3}
2 & 4 &\multicolumn{2}{c}{\;\;\;} \\
\cline{1-2}
\end{array}
\,.
\end{eqnarray}
Then tensors $T\in\calT_4 V$, $\hat{T}\in\calT_5 V$ fulfil
\begin{eqnarray}
T\in\calA(V) & \;\;\Leftrightarrow\;\; &
\textstyle{\frac{1}{12}}\,y_t^{\ast} T\;=\;T\,, \\
\hat{T}\in\calA'(V) & \;\;\Leftrightarrow\;\; &
\textstyle{\frac{1}{24}}\,y_{t'}^{\ast} \hat{T}\;=\;\hat{T}\,.
\end{eqnarray}
\end{Thm}
The group ring elements $\frac{1}{12}\,y_t^{\ast}$, $\frac{1}{24}\,y_{t'}^{\ast}$ are idempotents which are proportional to the essentially idempotent symmetrizers $y_t^{\ast}$, $y_{t'}^{\ast}$.
\vspace{0.5cm}

The group ring $\bbK[\calS_3]$ decomposes into the minimal 2-sided ideals ${\mathfrak a}_{(3)}, {\mathfrak a}_{(2\,1)}, {\mathfrak a}_{(1^3)}$. Whereas the 2-sided ideals ${\mathfrak a}_{(3)}, {\mathfrak a}_{(1^3)}\subset\bbK[\calS_3]$ have dimension 1 and define consequently unique symmetry classes of $\calT_3 V$, the 2-sided ideal ${\mathfrak a}_{(2\,1)}\subset\bbK[\calS_3]$ has dimension\footnote{The dimensions of the ${\mathfrak a}_{\lambda}$ can be calculated by means of the hook length formula.} 4 and contains an infinite set of minimal right ideals $\frakr$ which lead to an infinite set $\fS$ of symmetry classes $\calT_{\frakr}$ for tensors of $\calT_3 V$.
The set of generating idempotents for these right ideals $\frakr$ is infinite, too. In \cite{fie03b} we used {\itshape discrete Fourier transforms} to determine a family of primitive generating idempotents of the above minimal right ideals $\frakr\subset\bbK[\calS_3]$.

We denote by $\bbK^{d\times d}$ the set of all $d\times d$-matrices of elements of $\bbK$.
\begin{Def}
A {\it discrete Fourier transform}\footnote{See M. Clausen and U. Baum \cite{clausbaum1,clausbaum2} for details about fast discrete Fourier transforms.} for $\calS_r$ is an isomorphism
\begin{eqnarray}
D : \; \bbK [\calS_r] & \rightarrow &
\bigotimes_{\lambda \vdash r} {\bbK}^{d_{\lambda} \times d_{\lambda}} \label{e2.11}%
\end{eqnarray}
according to Wedderburn's theorem which maps the group ring $\bbK [\calS_r]$ onto an outer direct product $\bigotimes_{\lambda \vdash r} {\bbK}^{d_{\lambda} \times d_{\lambda}}$ of full matrix rings ${\bbK}^{d_{\lambda} \times d_{\lambda}}$. We denote by $D_{\lambda}$ the {\it natural projections}
$D_{\lambda} : \bbK [\calS_r] \rightarrow
{\bbK}^{d_{\lambda} \times d_{\lambda}}$.
\end{Def}
A discrete Fourier transform maps every group ring element $a\in\bbK[\calS_r]$
to a block diagonal matrix
\begin{eqnarray}
D :\;\;a\;=\;\sum_{p\in\calS_r}\,a(p)\,p & \mapsto &
\left(
\begin{array}{cccc}
A_{{\lambda}_1} & & & 0 \\
 & A_{{\lambda}_2} & & \\
 & & \ddots & \\
0 & & & A_{{\lambda}_k} \\
\end{array}
\right)\,.
\end{eqnarray}
The matrices $A_{\lambda}$ are numbered by the partitions $\lambda\vdash r$. The dimension $d_{\lambda}$ of every matrix $A_{\lambda}\in{\bbK}^{d_{\lambda} \times d_{\lambda}}$ can be calculated from the Young frame belonging to $\lambda\vdash r$ by means of the {\itshape hook length formula}. For $r = 3$ we have
\begin{center}
\begin{tabular}{c|ccc}
$\lambda$ & $(3)$ & $(2\,1)$ & $(1^3)$ \\
\hline
$d_{\lambda}$ & 1 & 2 & 1 \\
\end{tabular}
\end{center}
The inverse discrete Fourier transform is given by
\begin{Prop}\footnote{See M. Clausen and U. Baum \cite[p.81]{clausbaum1}.}
If $D : \bbK [\calS_r]\rightarrow
\bigotimes_{\lambda \vdash r} {\bbK}^{d_{\lambda} \times d_{\lambda}}$
is a discrete Fourier transform for $\bbK[\calS_r]$, then we have for every
$a\in\bbK[\calS_r]$
\begin{eqnarray}
\forall\,p\in\calS_r:\;\;\;a(p) & = & \frac{1}{r!}\,\sum_{\lambda\vdash r}\,
d_{\lambda}\,\mathrm{trace}\left\{D_{\lambda}(p^{-1})\cdot D_{\lambda}(a)\right\} \label{e2.13}\\
 & = & \frac{1}{r!}\,\sum_{\lambda\vdash r}\,
d_{\lambda}\,\mathrm{trace}\left\{D_{\lambda}(p^{-1})\cdot A_{\lambda}\right\}\,. \nonumber
\end{eqnarray}
\end{Prop}
In our considerations we are interested in the matrix ring $\bbK^{2\times 2}$ which corresponds to the $(2\,1)$-equivalence class of minimal right ideals $\frakr\subset\bbK[\calS_3]$. In \cite{fie03b} we proved
\begin{Prop}
Every minimal right ideal $\frakr\subset\bbK^{2\times 2}$ is generated by exactly one of the following (primitive) idempotents
\begin{eqnarray} \label{e2.14}%
Y\;:=\;\left(
\begin{array}{cc}
0 & 0 \\
0 & 1 \\
\end{array}
\right) 
 & \;\;\;or\;\;\; &
X_{\nu}\;:=\;\left(
\begin{array}{cc}
1 & 0 \\
\nu & 0 \\
\end{array}
\right)\;\;\;,\;\;\;\nu\in\bbK\,.
\end{eqnarray}
\end{Prop}
Using an inverse discrete Fourier transform we can
determine the primitive idempotents $\eta, {\xi}_{\nu}\in\bbK[\calS_3]$ which correspond to $Y$, $X_{\nu}$ in (\ref{e2.14}). We calculated these idempotents by means of the tool
\verb|InvFourierTransform| of the Mathematica package {\sf PERMS} \cite{fie10} (see also \cite[Appendix B]{fie16}.) This calculation can be carried out also by the program package {\sf SYMMETRICA} \cite{kerbkohn2,kerbkohnlas}.
\begin{Prop}
Let us use Young's natural representation\footnote{Three discrete Fourier transforms (\ref{e2.11}) are known for symmetric groups $\calS_r$: (1) {\itshape Young's natural representation}, (2) {\itshape Young's seminormal representation} and (3) {\itshape Young's orthogonal representation}. See \cite{boerner,boerner2,kerber,clausbaum1}. A short description of (1) and (2) can be found in \cite[Sec.I.5.2]{fie16}. All three discrete Fourier transforms are implemented in the program package {\sf SYMMETRICA} \cite{kerbkohn2,kerbkohnlas}. {\sf PERMS} \cite{fie10}  uses the natural representation. The fast {DFT}-algorithm of M. Clausen and U. Baum \cite{clausbaum1,clausbaum2} is based on the seminormal representation.} of $\calS_3$ as discrete Fourier transform. If we apply the Fourier inversion formula {\rm (\ref{e2.13})} to a $4\times 4$-block matrix
\begin{eqnarray}
\left(
\begin{array}{ccc}
A_{(3)} & & 0 \\
 & A_{(2\,1)} & \\
0 & & A_{(1^3)} \\
\end{array}
\right)
 & = &
\left(
\begin{array}{ccc}
0 & 0 & 0 \\
0 & A & 0 \\
0 & 0 & 0 \\
\end{array}
\right)
\end{eqnarray}
where $A$ is equal to $X_{\nu}$ or $Y$ in {\rm (\ref{e2.14})}, then we obtain the following idempotents of $\bbK[\calS_3]$:
\begin{eqnarray} \label{e2.16}%
X_{\nu}\;\;\;\Rightarrow\;\;\;{\xi}_{\nu} & := &
{\textstyle\frac{1}{3}}\,\{[1,2,3] + \nu [1,3,2] + (1-\nu)[2,1,3] \\
 & & - \nu [2,3,1] + (-1+\nu)[3,1,2] - [3,2,1]\} \nonumber \\
Y\;\;\;\Rightarrow\;\;\;\eta & := &
{\textstyle\frac{1}{3}}\,\{[1,2,3] - [2,1,3] - [2,3,1] + [3,2,1]\}\,. \label{e2.18}%
\end{eqnarray}
\end{Prop}
\begin{Rem}
It is interesting to clear up the connection of the idempotents ${\xi}_{\nu}$ and $\eta$ with Young symmetrizers.
A simple calculation shows that
\begin{eqnarray}
{\xi}_0\;=\;{\textstyle\frac{1}{3}}\,y_{t_1}
 & , &
\eta\;=\;{\textstyle\frac{1}{3}}\,y_{t_2}
\end{eqnarray}
where $y_{t_1}$ and $y_{t_2}$ are the Young symmetrizers of the tableaux
\begin{eqnarray*}
t_1\;=\;
\begin{array}{|c|c|c}
\cline{1-2}
1 & 2 & \\
\cline{1-2}
3 & \multicolumn{2}{c}{\;\;\;} \\
\cline{1-1}
\end{array}
 & , &
t_2\;=\,
\begin{array}{|c|c|c}
\cline{1-2}
1 & 3 & \\
\cline{1-2}
2 & \multicolumn{2}{c}{\;\;\;} \\
\cline{1-1}
\end{array}
\,.
\end{eqnarray*}
\end{Rem}
\begin{Rem} \label{rem2.12}%
In \cite[Thm 6.1]{fie03c} we showed that the symmetry classes $\calT_{\frakr_s}$, $\calT_{\frakr_a}$ of the tensors (\ref{e1.18}) correspond to the following values of the parameter $\nu$ in (\ref{e2.16}):
\begin{eqnarray}
{\xi}_{\nu}\cdot\bbK[\calS_3] = \frakr_s\;\Leftrightarrow\;\nu = 0
& \;,\; &
{\xi}_{\nu}\cdot\bbK[\calS_3] = \frakr_a\;\Leftrightarrow\;\nu = 2\,.
\end{eqnarray}
\end{Rem}
\vspace{10pt}

\section{Proof of the Theorems \ref{thm1.4}, \ref{thm1.5} and of Lemma \ref{lem1.6}}
\label{sec3}%
\subsection{Proof of Theorem \ref{thm1.4}}
In \cite{fie03b} the statements of Theorem \ref{thm1.4} about the products (a'), (b'), (c') were proved by an application of the \textit{Littlewood-Richardson rule}\footnote{See the
references \cite{kerber,kerber3,jameskerb,littlew1,mcdonald,full4,fultharr} for the Littlewood-Richardson rule. See also \cite{fie16}.}. Using \cite[Prop.2.10]{fie03b} and the arguments of \cite[Sec.3.1]{fie03b} we can prove the assertion about the products (a) of Theorem \ref{thm1.4} in the same way.

For tensor products $U\otimes w$ and $w\otimes U$ we have to investigate \textit{Littlewood-Richardson products} $[\lambda][1]\sim [\lambda]\,\#\,[1]\uparrow\calS_4$, $\lambda\vdash 3$. The three partitions
$(3), (2\,1), (1^3)\vdash 3$ of $3$ lead to the three Littlewood-Richardson products
\begin{eqnarray*}
[3][1] & \sim & [4] + [3\,1]\,, \\
\,[2\,1][1] & \sim & [3\,1] + [2^2] + [2\,1^2]\,, \\
\,[1^3][1] & \sim & [2\,1^2] + [1^4]\,.
\end{eqnarray*}
But then we obtain by means of the arguments of \cite[Sec.3.1]{fie03b}, that the symmetry class of $U$ must belong to the partition $(2\,1)\vdash 3$ since only the product $[2\,1][1]$ contains a part $[2^2]$. The part $[2^2]$ describes a minimal left ideal $\frakl\subset\bbK[\calS_4]$ which lies in the same equivalence class of minimal left ideals as the left ideal $\bbK[\calS_4]\cdot y_t$ generated by the Young symmetrizer $y_t$ from Theorem \ref{thm2.6} with the Young frame $(2^2)\vdash 4$.
Consequently, only for tensors $T = U\otimes w$ with a part $[2^2]$ a relation $0\not= T_b\cdot y_t = (y_t^{\ast}T)_b$, $b\in V^4$, or equivalently $y_t^{\ast}(U\otimes w)\not= 0$ is possible. $\Box$

\subsection{Proof of Theorem \ref{thm1.5}}
The equivalence of the Statements (\ref{i1.5.2}), (\ref{i1.5.3}), (\ref{i1.5.4}) of Theorem \ref{thm1.5} follows from \cite[Thm.1.10]{fie03b}. The equivalence of (\ref{i1.5.1}) and (\ref{i1.5.4}) can also be proved by means of conclusions which were used in the proof of \cite[Thm.1.10]{fie03b}.

For a treatment of expressions $y_{t}^{\ast}(U\otimes w)$ we form the following elements of $\bbK[\calS_4]$:
\begin{eqnarray}
{\sigma}_{\nu} & := & y_{t}^{\ast}\cdot {\xi}_{\nu}'
\;\;\;,\;\;\;
{\rho}\;:=\;y_{t}^{\ast}\cdot {\eta}'\\
{\xi}_{\nu} & \mapsto & {\xi}_{\nu}'\in\bbK[\calS_4]
\;\;\;,\;\;\;
{\eta}\;\mapsto\;{\eta}'\in\bbK[\calS_4]\,. \label{e3.2}
\end{eqnarray}

Formula (\ref{e3.2}) denotes the embedding of ${\xi}_{\nu}\,,\,\eta\in\bbK[\calS_3]$ into $\bbK[\calS_4]$ which is induced by the mapping
$\calS_3\rightarrow\calS_4\;,\;
[i_1,i_2,i_3]\mapsto [i_1,i_2,i_3,4]$.

Using the Mathematica package {\sf PERMS} \cite{fie10} we verified that
\begin{center}
\fbox{
${\rho}\not= 0 \;\;\;\;\;\;\mathrm{and}\;\;\;\;\;\;
{\sigma}_{\nu}\not= 0\;\Leftrightarrow\;\nu\not=\frac{1}{2}\,.
$}
\end{center}
If ${\rho}\not= 0$, ${\sigma}_{\nu}\not= 0$ then the
minimal right ideals $y_{t}^{\ast}\cdot\bbK[\calS_4]$, ${\rho}\cdot\bbK[\calS_4]$ and  ${\sigma}_{\nu}\cdot\bbK[\calS_4]$
coincide, i.e. the symmetry operators ${\rho}$, ${\sigma}_{\nu}$ can be used to define the symmetry class $\calA(V)$ of algebraic curvature tensors.
A tensor $T\in\calT_4 V$ is an algebraic curvature tensor iff a tensor $T'\in\calT_4 V$ exists such that
$T = {\rho} T'$ or $T = {\sigma}_{\nu} T'$ (if $\nu\not=\frac{1}{2}$).

Now we represent the tensor $T'$ by by a finite sum of product tensors
\begin{eqnarray}
T' & = & \sum_{(M,w)\in\calP}\,M\otimes w\;\;\;,\;\;\;
\calP\subset\calT_3 V\times\calT_1 V\;\mathrm{finite}.
\end{eqnarray}
The symmetrizations
${\xi}_{\nu}'(M\otimes w)$,
${\eta}'(M\otimes w)$, 
lead to product tensors
$U\otimes w$, where $U$ has a $(2\,1)$-symmetry, defined by ${\xi}_{\nu}$ or $\eta$.

Thus we can write every tensor
${\xi}_{\nu}'T'$ and
${\eta}'T'$ 
as a finite sum of suitable tensors 
$U\otimes w$.
An application of $y_{t}^{\ast}$ shows that every algebraic curvature tensor $T$ can be expressed by a finite sum of tensors $y_{t}^{\ast}(U\otimes w)$ if $\nu\not=\frac{1}{2}$.
Finally we can verify by a simple calculation that the idempotent ${\xi}_{\nu}$ is equal to the idempotent (\ref{e1.14}) if $\nu = \frac{1}{2}$. This finishes the proof of Theorem \ref{thm1.5}. \textsf{Mathematica} notebooks of \textsf{PERMS}-calculations for this proof can be found on the web page \cite{fie21}. $\Box$

\subsection{Proof of Lemma \ref{lem1.6}}
The formulas (\ref{e1.16}) and (\ref{e1.17}) in Lemma \ref{lem1.6} were proved in \cite[Appendix A]{fie03c} by computer calculations using the \textsf{Mathematica} packages \textsf{Ricci} \cite{ricci3} and \textsf{PERMS} \cite{fie10}. A proof of (\ref{e1.15}) can be obtained in the same way. The calculation stored in the notebook \cite[part1.nb]{fie21} yield (\ref{e1.15}) and the following expression of 16 terms for the coordinates of $\frac{1}{12}\,y_t^{\ast}(U\otimes w)$:
{\small
\begin{eqnarray}
 & &
\frac{1}{12} \,\tensor{U}{\down{j}\down{k}\down{l}} 
    \,\tensor{w}{\down{i}} - 
  \frac{1}{12} \,\tensor{U}{\down{j}\down{l}\down{k}} 
    \,\tensor{w}{\down{i}} - 
  \frac{1}{12} \,\tensor{U}{\down{k}\down{l}\down{j}} 
    \,\tensor{w}{\down{i}} + 
  \frac{1}{12} \,\tensor{U}{\down{l}\down{k}\down{j}} 
    \,\tensor{w}{\down{i}} - 
  \frac{1}{12} \,\tensor{U}{\down{i}\down{k}\down{l}} 
    \,\tensor{w}{\down{j}} + 
  \frac{1}{12} \,\tensor{U}{\down{i}\down{l}\down{k}} 
    \,\tensor{w}{\down{j}} + \nonumber \\ & & 
  \frac{1}{12} \,\tensor{U}{\down{k}\down{l}\down{i}} 
    \,\tensor{w}{\down{j}} - 
  \frac{1}{12} \,\tensor{U}{\down{l}\down{k}\down{i}} 
    \,\tensor{w}{\down{j}} - 
  \frac{1}{12} \,\tensor{U}{\down{i}\down{j}\down{l}} 
    \,\tensor{w}{\down{k}} + 
  \frac{1}{12} \,\tensor{U}{\down{j}\down{i}\down{l}} 
    \,\tensor{w}{\down{k}} + 
  \frac{1}{12} \,\tensor{U}{\down{l}\down{i}\down{j}} 
    \,\tensor{w}{\down{k}} - 
  \frac{1}{12} \,\tensor{U}{\down{l}\down{j}\down{i}} 
    \,\tensor{w}{\down{k}} + \label{e3.4}\\ & &
  \frac{1}{12} \,\tensor{U}{\down{i}\down{j}\down{k}} 
    \,\tensor{w}{\down{l}} - 
  \frac{1}{12} \,\tensor{U}{\down{j}\down{i}\down{k}} 
    \,\tensor{w}{\down{l}} - 
  \frac{1}{12} \,\tensor{U}{\down{k}\down{i}\down{j}} 
    \,\tensor{w}{\down{l}} + 
  \frac{1}{12} \,\tensor{U}{\down{k}\down{j}\down{i}} 
    \,\tensor{w}{\down{l}} \nonumber
\end{eqnarray}
}
\vspace{10pt}

\section{Short formulas for algebraic curvature tensors $\fR$}
\label{sec4}%
\subsection{The reduction method}
In this section we begin to construct short coordinate representations of tensors $y_t^{\ast}(U\otimes w)$ considered in Theorem \ref{thm1.5} and Lemma \ref{lem1.6}.
Formula (\ref{e3.4}) represents the coordinates of $\frac{1}{12}\,y_t^{\ast}(U\otimes w)$ by a relatively long polynomial
\begin{eqnarray} \label{e4.1}%
\mathfrak{P}_{i_1\ldots i_4} & := & {\textstyle\frac{1}{12}}\,y_t^{\ast}(U\otimes w)_{i_1\ldots i_4}\;=\;
\sum_{p\in\calS_4}\,c_p\,U_{i_{p(1)} i_{p(2)} i_{p(3)}} w_{i_{p(4)}}\;\;\;,\;\;\;c_p\in\bbK
\end{eqnarray}
in the coordinates of $U$ and $w$. In \cite[Sec.3,4]{fie03c} we developed a method to reduce the length of the coordinate representation of tensors $y_{t'}(U\otimes S)$ and $y_{t'}(U\otimes A)$ from Theorem \ref{thm1.5}. We can also use this method for a reduction of the length of (\ref{e4.1}). Here is a summary of the method from \cite[Sec.3,4]{fie03c}.

A central role plays the following
\begin{Prop}
Let $\frakr\subset\bbK[\calS_r]$ be a $d$-dimensional right ideal that defines a symmetry class $\calT_{\frakr}$ of tensors $T\in\calT_r V$.
If a basis
$\calB = \{ h_1 , \ldots , h_d \}$ of the left ideal $\frakl = \frakr^{\ast}$ is known, then
the coefficients $x_p\in\bbK$ for linear identities
\begin{eqnarray}
\sum_{p\in\calS_r}\,x_p\,T_{i_{p(1)}\ldots i_{p(r)}} & = & 0
\end{eqnarray}
satisfied by the coordinates of the $T\in\calT_{\frakr}$
can be obtained
from the linear $(d\times r!)$-equation system
\begin{eqnarray} \label{e4.3}%
\sum_{p \in {\calS}_r} h_i(p)\,x_p & = & 0 \hspace{1cm}(i = 1 , \ldots , d)\,.
\end{eqnarray}
\end{Prop}
For the tensors $U$ considered in Theorem \ref{thm1.5} the system (\ref{e4.3}) has a $(2\times 6)$-coefficient matrix since $d = \dim\frakr = 2$ for the right ideal $\frakr$ defining the symmetry class of $U$.
Using discrete Fourier transforms and results from \cite{fie16} we proved in \cite{fie03c}
\begin{Lem}
The left ideals $\bbK[\calS_3]\cdot {\xi}_{\nu}^{\ast}$, $\bbK[\calS_3]\cdot {\eta}^{\ast}$, given by the idempotents {\rm (\ref{e2.16})}, {\rm (\ref{e2.18})}, possess bases which lead to the following coefficient matrices in {\rm (\ref{e4.3})}
\begin{eqnarray}
\hspace*{1cm} & {\xi}_{\nu}\;\Rightarrow &
\frac{1}{9}\,\left(
\begin{array}{rrrrrr}
4 - 2\nu & -2 + 4\nu & 4 - 2\nu & -2 + 4\nu & -2 - 2\nu & -2 - 2\nu \\
-2 + 4\nu & 4 - 2\nu & -2 - 2\nu & -2 - 2\nu & 4 - 2\nu & -2 + 4\nu \\
\end{array}
\right) \label{e4.4}\\
\hspace*{1cm} & \eta\;\Rightarrow &
\frac{1}{9}\,\left(
\begin{array}{rrrrrr}
-1 & 2 & -1 & 2 & -1 & -1 \\
2 & -1 & -1 & -1 & -1 & 2 \\
\end{array}
\right)\,. \label{e4.5}%
\end{eqnarray}
Here $\nu\in\bbK$ is arbitrary. Further we use the following correspondence $a\leftrightarrow p_a$ between the column index $a$ in {\rm (\ref{e4.4})}, {\rm (\ref{e4.5})} and permutations $p_a\in\calS_3$:
\begin{eqnarray} \label{e4.6}%
 & &
\begin{array}{|c|c|c|c|c|c|c|}
\hline
a & 1 & 2 & 3 & 4 & 5 & 6 \\
\hline
p_a & [1,2,3] & [1,3,2] & [2,1,3] & [2,3,1] & [3,1,2] & [3,2,1] \\
\hline
\end{array}\,.
\end{eqnarray}
\end{Lem}

For $a, b\in\{1,\ldots , 6\}$, $a < b$, we denote by $\calP_{ab}$ the 2-set
$\calP_{ab} := \{p_a , p_b \}$ of permutations from $\calS_3$ which correspond to $a , b$ via (\ref{e4.6}). Furthermore we write $\Delta_{\calP_{ab}}$ for the determinant of the $(2\times 2)$-submatrix of (\ref{e4.4}) or (\ref{e4.5}) whose columns correspond to $p_a , p_b$.
\begin{Proc} \label{proc4.3}%
Consider a symmetry class $\calT_{\frakr}$ of tensors $U\in\calT_3 V$ defined by ${\xi}_{\nu}$ or $\eta$, and the corresponding equation system (\ref{e4.3}) with coefficient matrix (\ref{e4.4}) or (\ref{e4.5}). Then carry out the following steps for every set $\calP_{ab}$:
\begin{enumerate}
\item{Check the condition $\Delta_{\calP_{ab}}\not= 0$. If $\Delta_{\calP_{ab}} = 0$, then skip the steps (\ref{i4.3.2}), \ldots , (\ref{i4.3.4}) for $\calP_{ab}$.}
\item{If $\Delta_{\calP_{ab}}\not= 0$, then, for every $\bar{p}\in\calS_3\setminus\calP_{ab}$, determine the solution $x_p^{(\bar{p})}$ of (\ref{e4.3}) which fulfils $x_{\bar{p}}^{(\bar{p})} = 1$ and $x_{p}^{(\bar{p})} = 0$ for all $p\in\calS_3\setminus (\calP_{ab}\cup\{\bar{p}\})$. \label{i4.3.2}}
\item{Use the $x_p^{(\bar{p})}$ of step (\ref{i4.3.2}) to form identities
\begin{eqnarray} \label{e4.7}%
0 & = & \sum_{p\in\calP_{ab}}\,x_p^{(\bar{p})}\,U_{i_{p(1)} i_{p(2)} i_{p(3)}} +
U_{i_{\bar{p}(1)} i_{\bar{p}(2)} i_{\bar{p}(3)}}\hspace{1cm}(\bar{p}\in{\calS_3}\setminus\calP_{ab})\,.
\end{eqnarray}
}
\item{Interpret $\{i_{p(1)}, i_{p(2)}, i_{p(3)}\}$ as a permuted arrangement of a lexicographically ordered sequence $\{i_1, i_2, i_3\}$ of index names. Use (\ref{e4.7}) to express all coordinates $U_{i_{\bar{p}(1)} i_{\bar{p}(2)} i_{\bar{p}(3)}}\,,\,\bar{p}\in{\calS_3}\setminus\calP_{ab}$, by the coordinates $U_{i_{p(1)} i_{p(2)} i_{p(3)}}\,,\,p\in\calP_{ab}$, in (\ref{e3.4}). \label{i4.3.4}}
\end{enumerate}
\end{Proc}
For instance, let us consider (\ref{e4.4}). Then the set
$\calP_{1\,2} = \{[1,2,3]\,,\,[1,3,2]\}$ leads to the determinant
${\Delta}_{\calP_{1\,2}}(\nu) = {\textstyle\frac{4}{27}}\,(1 - \nu)(1 + \nu)$
which has the roots ${\nu}_1 = 1$ and ${\nu}_2 = -1$. For $\nu\not\in\{1\,,\,-1\}$ we obtain the identities
\begin{eqnarray} \label{e4.8}%
 & &
\begin{array}{ccccccccccc}
- & \frac{{\nu}^2 - \nu + 1}{{\nu}^2 - 1}\,U_{ijk} & + & \frac{2\nu - 1}{{\nu}^2 - 1}\,U_{ikj}        & +  &         &         &         & U_{kji} & = & 0 \\
  & \frac{{\nu}^2 - 2\nu}{{\nu}^2 - 1}\,U_{ijk} & + & \frac{{\nu}^2 - \nu + 1}{{\nu}^2 - 1}\,U_{ikj} & + &         &         & U_{kij} &         & = & 0 \\
  & \frac{2\nu - 1}{{\nu}^2 - 1}\,U_{ijk}        & - & \frac{{\nu}^2 - \nu + 1}{{\nu}^2 - 1}\,U_{ikj} & + &         & U_{jki} &         &         & = & 0 \\
  & \frac{{\nu}^2 - \nu + 1}{{\nu}^2 - 1}\,U_{ijk} & + & \frac{{\nu}^2 - 2\nu}{{\nu}^2 - 1}\,U_{ikj} & + & U_{jik} &         &         &         & = & 0 \\
\end{array}\;.
\end{eqnarray}
There exist 15 subsets $\calP_{ab} = \{p_a,p_b\}\subset\calS_3$ and consequently 15 systems of type (\ref{e4.8}) for $U$ with respect to (\ref{e4.4}). The matrix (\ref{e4.5}) leads to 12 systems (\ref{e4.8}) because $\Delta_{\calP_{1\,6}} = \Delta_{\calP_{2\,4}} = \Delta_{\calP_{3\,5}} = 0$ (see \cite[Table 1]{fie03c}).

\subsection{Proof of Theorem \ref{thm1.8}, (1a)}
We carried out Procedure \ref{proc4.3} in computer calculations using the \textsf{Mathematica} packages \textsf{Ricci} \cite{ricci3} and \textsf{PERMS} \cite{fie10}. The \textsf{Mathematica} notebooks are available on the web page \cite{fie21}. In all cases with $\Delta_{\calP_{ab}}\not= 0$ and $\nu\not= \frac{1}{2}$ we obtained a reduction of (\ref{e3.4}) to 8 terms both for (\ref{e4.4}) and for (\ref{e4.5}). The roots of $\Delta_{\calP_{ab}}(\nu)$ for (\ref{e4.4}) are given in \cite[Table 2]{fie03c}. $\Box$

\subsection{Proof of Theorem \ref{thm1.8}, (1b) and (2)}
For symmetry classes of tensors $U$ described by ${\xi}_{\nu}$ and (\ref{e4.4}) a further reduction of the length of (\ref{e3.4}) is possible.
When we use a system of linear identities (\ref{e4.7}) to reduce the length of (\ref{e4.1}), (\ref{e3.4}) then we obtain a sum with a structure
\begin{eqnarray} \label{e4.9}%
\mathfrak{P}_{i_1\ldots i_4}^{\rm red} & = & \sum_{q\in\calS_4}\,\frac{P_q^{\calP_{ab}}(\nu)}{Q_q^{\calP_{ab}}(\nu)}\,U_{i_{q(1)} i_{q(2)} i_{q(3)}} w_{i_{q(4)}}\,,
\end{eqnarray}
where $P_q^{\calP_{ab}}(\nu)$ and $Q_q^{\calP_{ab}}(\nu)$ are polynomials, because the entries of (\ref{e4.4}) are polynomials. For instance, the identities (\ref{e4.8}) belonging to
$\calP_{1\,2}$ transform (\ref{e3.4}) into
{\small
\begin{eqnarray}
 {\textstyle\frac{1}{12}}\,y_t^{\ast}(U\otimes w)_{ijkl} & = &
\hspace{-0.2cm}{\textstyle\frac{-1 + 2\,\nu}{12\,\left( -1 + \nu \right) }} 
   \,\tensor{U}{\down{j}\down{k}\down{l}} \,\tensor{w}{\down{i}} - 
  {\textstyle\frac{-1 + 2\,\nu}{12\,\left( -1 + \nu \right) }} 
   \,\tensor{U}{\down{j}\down{l}\down{k}} \,\tensor{w}{\down{i}} - 
  {\textstyle\frac{-1 + 2\,\nu}{12\,\left( -1 + \nu \right) }} 
   \,\tensor{U}{\down{i}\down{k}\down{l}} \,\tensor{w}{\down{j}} + \nonumber \\ & &
\hspace{-0.2cm}{\textstyle\frac{-1 + 4\,\nu - 4\,{\nu}^2}{12\,\left( -1 + \nu \right) \,\left( 1 + \nu \right) }} 
   \,\tensor{U}{\down{i}\down{j}\down{l}} \,\tensor{w}{\down{k}} + 
  {\textstyle\frac{-2 + 5\,\nu - 2\,{\nu}^2}{12\,\left( -1 + \nu \right) \,\left( 1 + \nu \right) }} 
   \,\tensor{U}{\down{i}\down{l}\down{j}} \,\tensor{w}{\down{k}} + 
  {\textstyle\frac{-1 + 2\,\nu}{12\,\left( -1 + \nu \right) }} 
   \,\tensor{U}{\down{i}\down{l}\down{k}} \,\tensor{w}{\down{j}} + \nonumber \\ & &
\hspace{-0.2cm}{\textstyle\frac{1 - 4\,\nu + 4\,{\nu}^2}{12\,\left( -1 + \nu \right) \,\left( 1 + \nu \right) }} 
   \,\tensor{U}{\down{i}\down{j}\down{k}} \,\tensor{w}{\down{l}} + 
  {\textstyle\frac{2 - 5\,\nu + 2\,{\nu}^2}{12\,\left( -1 + \nu \right) \,\left( 1 + \nu \right) }} 
   \,\tensor{U}{\down{i}\down{k}\down{j}} \,\tensor{w}{\down{l}}\,. \label{e4.10}%
\end{eqnarray}
}

Now, for every $\calP_{ab}$ we determine the set $N_{\calP_{ab}}$ of all roots $\nu$ of the $P_q^{\calP_{ab}}(\nu)$ in (\ref{e4.9}) which are different from the roots of ${\Delta}_{\calP_{ab}}(\nu)$, $Q_q^{\calP_{ab}}(\nu)$ and from $\nu = \frac{1}{2}$.
If we set such a $\nu\in N_{\calP_{ab}}$ into (\ref{e4.9}), the length of (\ref{e4.9}) will decrease. For example, (\ref{e4.10}) yield $N_{\calP_{1\,2}} = \{2\}$. The root $\nu = 2$ reduces (\ref{e4.10}) to the 6 terms
{\small
\begin{eqnarray}
 & &
{\textstyle\frac{1}{4}} \,\tensor{U}{\down{j}\down{k}\down{l}} 
    \,\tensor{w}{\down{i}} - 
  {\textstyle\frac{1}{4}} \,\tensor{U}{\down{j}\down{l}\down{k}} 
    \,\tensor{w}{\down{i}} - 
  {\textstyle\frac{1}{4}} \,\tensor{U}{\down{i}\down{k}\down{l}} 
    \,\tensor{w}{\down{j}} + 
  {\textstyle\frac{1}{4}} \,\tensor{U}{\down{i}\down{l}\down{k}} 
    \,\tensor{w}{\down{j}} - 
  {\textstyle\frac{1}{4}} \,\tensor{U}{\down{i}\down{j}\down{l}} 
    \,\tensor{w}{\down{k}} + 
  {\textstyle\frac{1}{4}} \,\tensor{U}{\down{i}\down{j}\down{k}} 
    \,\tensor{w}{\down{l}}\,.
\end{eqnarray}
}
Using \textsf{Ricci} \cite{ricci3} and \textsf{PERMS} \cite{fie10}, we determined all sets $N_{\calP_{ab}}$ and all resulting length reductions of (\ref{e4.9}). Table 1 shows the results. The minimal length of $\mathfrak{P}_{i_1\ldots i_4}^{\rm red}$ which we found is equal to 4. For example, the set $\calP_{1\,6}$ lead to $N_{\calP_{1\,6}} = \{-1\,,\,2\}$ and these two roots reduce $\frac{1}{12} y_t^{\ast}(U\otimes w)_{ijkl}$ to the following 4 expressions
{\small
\begin{eqnarray*}
\nu = -1 & \Rightarrow &
{\textstyle\frac{1}{4}} \,\tensor{U}{\down{j}\down{k}\down{l}} 
    \,\tensor{w}{\down{i}} - 
  {\textstyle\frac{1}{4}} \,\tensor{U}{\down{i}\down{k}\down{l}} 
    \,\tensor{w}{\down{j}} - 
  {\textstyle\frac{1}{4}} \,\tensor{U}{\down{l}\down{j}\down{i}} 
    \,\tensor{w}{\down{k}} + 
  {\textstyle\frac{1}{4}} \,\tensor{U}{\down{k}\down{j}\down{i}} 
    \,\tensor{w}{\down{l}}\\
\nu = \hspace{0.3cm}2 & \Rightarrow &
{\textstyle\frac{1}{4}} \,\tensor{U}{\down{l}\down{k}\down{j}} 
    \,\tensor{w}{\down{i}} - 
  {\textstyle\frac{1}{4}} \,\tensor{U}{\down{l}\down{k}\down{i}} 
    \,\tensor{w}{\down{j}} - 
  {\textstyle\frac{1}{4}} \,\tensor{U}{\down{i}\down{j}\down{l}} 
    \,\tensor{w}{\down{k}} + 
  {\textstyle\frac{1}{4}} \,\tensor{U}{\down{i}\down{j}\down{k}} 
    \,\tensor{w}{\down{l}}
\end{eqnarray*}
}
Similar tables for $y_{t'}^{\ast}(U\otimes S)$ and $y_{t'}^{\ast}(U\otimes A)$ are given in \cite{fie03c}. They yield the minimal lengths 12 and 10 in Statement (1b) of Theorem \ref{thm1.8}.
\begin{table}[t]
\begin{center}
\begin{tabular}{|c|c|c|}
\hline
$\calP_{ab}$ & $N_{\calP_{ab}}$ & length $\mathfrak{P}_{i_1\ldots i_4}^{\rm red}$ \\
\hline
 12  &    2  &      6 \\
\hline
 13  &    -1 &      6 \\
\hline
 14  &    -1 &      6 \\
     &    2  &      4 \\
\hline
 15  &    -1 &      4 \\
     &    2  &      6 \\
\hline
 16  &    -1 &      4 \\
     &    2  &      4 \\
\hline
 23  &    -1 &      6 \\
     &    2  &      6 \\
\hline
 24  &    -1 &      6 \\
     &    2  &      6 \\
\hline
\end{tabular}
\hspace{40pt}
\begin{tabular}{|c|c|c|}
\hline
$\calP_{ab}$ & $N_{\calP_{ab}}$ & length $\mathfrak{P}_{i_1\ldots i_4}^{\rm red}$ \\
\hline
 25  &    -1 &      4 \\
\hline
 26  &    -1 &      4  \\
     &    2  &      6  \\
\hline
 34  &    2  &      4  \\
\hline
 35  &    -1 &      6  \\
     &    2  &      6  \\
\hline
 36  &    -1 &      6  \\
     &    2  &      4  \\
\hline
 45  &    -1 &      6  \\
     &    2  &      6  \\
\hline
 46  &    -1 &      6  \\
\hline
 56  &    2  &      6  \\
\hline
\end{tabular}
\vspace{3mm}
\caption{The lengths of $\mathfrak{P}_{i_1\ldots i_4}^{\rm red}$ for an $U$ from a ${\xi}_{\nu}$-symmetry class, where $\nu$ is an allowed root of a $P_q^{\calP_{ab}}(\nu)$.}
\end{center}
\end{table}
Furthermore Table 1 and the tables in \cite{fie03c} show that Statement (2) of Theorem \ref{thm1.8} is valid exactly for the two right ideals
$\frakr = {\xi}_{\nu}\cdot\bbK[\calS_3]$ which belong to $\nu = -1$ and $\nu = 2$. $\Box$

\subsection{Proof of Theorem \ref{thm1.8}, (3) and Theorem \ref{thm1.10}}
In \cite[Sec.5]{fie03c} we proved that 
a symmetry class $\calT_{\frakr}$ of a minimal right ideal $\frakr\subset\bbK[\calS_3]$ admits an index commutation symmetry of the tensors $U\in\calT_{\frakr}$ iff $\frakr = \eta\cdot\bbK[\calS_3]$ or  $\frakr = {\xi}_{\nu}\cdot\bbK[\calS_3]$ with $\nu\in\{0, 1, 2, \frac{1}{2}, -1, e^{\imath\pi/3}, e^{-\imath\pi/3}\}$. The minimal lengths of Theorem \ref{thm1.8}, (1b) occured in following cases:
\begin{eqnarray*}
\mathrm{for}\;y_t^{\ast}(U\otimes w)\,,\,y_{t'}^{\ast}(U\otimes A)
 & \Leftrightarrow & \frakr = {\xi}_{\nu}\cdot\bbK[\calS_3]\;\mathrm{with}\;\nu\in\{-1\,,\,2\}\,, \\
\mathrm{for}\;y_{t'}^{\ast}(U\otimes S)
 & \Leftrightarrow & \frakr = \eta\cdot\bbK[\calS_3]\;\mathrm{or}\;\frakr = {\xi}_{\nu}\cdot\bbK[\calS_3]\,,\,\nu\in\{-1\,,\,0\,,\,1\,,\,2\}\,.
\end{eqnarray*}
This proves Theorem \ref{thm1.8}, (3). Furthermore, Theorem \ref{thm1.10} follows from Remark \ref{rem2.12}. $\Box$


\end{document}